\documentclass[preprint,12pt]{elsarticle}

\usepackage{amssymb}
\usepackage{amsthm}
\usepackage{graphics}
\usepackage{epsfig}
\usepackage{amssymb}
\usepackage{hyperref}
\usepackage{graphicx}
\usepackage{subfigure}
\usepackage{color}
\usepackage{tikz}
\usetikzlibrary{patterns}
\usepackage{hyperref}
\usepackage{a4wide}
\usepackage{amsmath}
\usepackage{amsfonts}
\usepackage{enumerate}
\newcommand{\pf}{\noindent {\bf Proof: }}

\newtheorem*{theoremaux}{Theorem \theoremauxnum}
\gdef\theoremauxnum{1}

\newtheorem{lemma}{\bf Lemma}[section]

\newtheorem{theorem}{\bf Theorem}[section]

\newtheorem{note}{\bf Note}[section]

\journal{~}

\begin{document}
	
	\begin{frontmatter}
		
		


		\title{On Automorphism Group of a Family of Symmetric Graphs}

		\author{Sucharita Biswas}
		\ead{biswas.sucharita56@gmail.com}

		\address{Department of Mathematics, Presidency University, Kolkata, India\\
			86/1, College Street, Kolkata, India.
		}
		
		\begin{abstract}
			
			In this paper, we study the automorphism group and geodesic transitivity of a family of vertex-transitive graphs $H(n,k)$, introduced by Fu-Tao Hu \textit{et.al.} in 2010. In the process, we address some naturally arising, unanswered questions from that paper.
		\end{abstract}
		
	\end{frontmatter}
	
	\section{Introduction}
	Let $G=(V, E)$ be a graph with vertex set $V$ and edge set $E$. An automorphism of $G$ is a permutation $\varphi$ on the vertex set $V$ such that for any two vertices $u,v$ we have $u \sim v$ if and only if $\varphi(u) \sim \varphi(v)$. Set of all automorphisms of a graph $G$ form a group under mapping composition. If the  automorphism group of $G$  acts transitively on the set of all vertices(edges)  then $G$ is called \textit{vertex(edge) transitive} graph. A graph $G$ is called \textit{symmetric} or \textit{arc transitive} if the automorphism group of $G$ acts transitively on the set of all arcs.
	Graphs and their symmetries play a very important role in modern algebraic graph theory. Starting from the 19th century, the automorphism group of different families of graphs were studied, like the automorphism group of generalized Petersen graphs \cite{Petersen}, hypercube \cite{hypercube}, bi-Cayley graphs \cite{bi-Cayley}, Johnson graph \cite{Johnson}, bipartite Kneser graph \cite{Kneser}, doubled Grassmann graphs \cite{Grassmann} and many others like \cite{Half-arc}, generalized Pappus graphs \cite{Pappus}, generalized Andr{\' a}sfai graphs \cite{Andrasfai},  \cite{Orthogonality}, symmetric cubic graphs \cite{Cubic} etc.  One can refer to the expository article \cite{survey} by P.J. Cameron for a quick reference on the automorphisms of finite graphs.
	
	A shortest path between any two vertices of a graph is called the \textit{geodesic} between those vertices. A geodesic is a \textit{s-geodesic} if the distance between those vertices is $s$. A graph $G$ is called \textit{locally s-geodesic transitive} if for each vertex $u$, the stabilizer of $u$  acts transitively on the set of all $i$-geodesics starting from $u$, for all $1 \leq i \leq s$. If the graph is vertex-transitive then $G$ is \textit{s-geodesic transitive}. Let $d(G)$ denote the diameter of $\Gamma$. If $s=d(G)$ then $G$ is called \textit{geodesic transitive}. The study of the geodesic transitivity of a graph is very useful to understand the symmetry of that graph. There are plenty of research works done on geodesic transitivity of the graphs \cite{thesis},\cite{jin},\cite{jin2} etc.
	
	Fu-Tao Hu \textit{et.al.} defined a new class of transitive Cayley graph $H(n,k)$ \cite{Hu_Wang_Xu} and studied some algebraic and topological properties.
	In \cite{Hu_Wang_Xu}, the authors mentioned the automorphism group for the case $(n,k)=(n,1)$. In this paper, we deal with the general case and compute the full automorphism group of $H(n,k)$ except for the case $n=2k$. In Section 2, we provide some necessary definitions and preliminaries, in Section 3  we deal with the automorphism group of $H(n,k)$ where $n=2k-1$ and in Sections 4 and 5 we deal with the case $n=2k+1$. Hereafter  Section 6 take care of the rest of the cases except $n=2k$. In the last section of this paper, we discuss the geodesic transitivity of $H(n,k)$ for all $n$ and $k$. 
	
	\section{Definitions and Preliminaries}
	
	Going as per the notations used in \cite{Hu_Wang_Xu}, $\Omega_n$ denotes the power set of the set $[n]=\{ 1,2,\ldots ,n \}$, $\Omega'_n$, $\Omega''_n$ and $\Omega_n ^i$ contains all subsets of $[n]$ with odd, even and `$i$' size respectively, where $i \in [n]$.  Let $X\in \Omega_n$, we denote complement  of $X$  in $[n]$ by $X^c$, i.e., $X^c=[n] \setminus X$.
	$H(n,k)$ is the graph with vertex set $\Omega_n$ and any two sets $X$, $Y \in \Omega_n$, $X \sim Y$ if and only if $|X \triangle Y|=k$, where $\triangle$ denotes the symmetric difference of sets. Therefore $H(n,0)$ is an empty graph with no edges and $H(n,n)=2^{n-1}K_2$, where $K_2$ is the complete graph with two vertices. Consider $H'(n,k)$ and $H''(n,k)$ be the subgraphs of $H(n,k)$ induced by $\Omega'$ and $\Omega''$ respectively.

	Let $\sigma \in S_n$ and $X=\{x_1,x_2, \ldots ,x_m\}$. Define 
	$\sigma(X)=\emptyset$ if $X=\emptyset$ and $\sigma(X)=\{\sigma(x_1),\sigma(x_2), \ldots ,\sigma(x_m)\}$ if $X \neq \emptyset$. 
	Let $X$ be an arbitrary fixed element in $\Omega_n$. Define $\rho_X: \Omega_n \rightarrow \Omega_n$ such that $\rho_X(Y)=X \triangle Y$, for all $Y \in \Omega_n$ and  $H_n=\{\rho_X : X \in \Omega_n\}$.

	{\theorem \label{H(n,n-k)}\cite{Hu_Wang_Xu} If $n$ is even and $k$ is odd, then $H(n,k) \cong H(n,n-k)$.}
	
	{\theorem \label{even k}\cite{Hu_Wang_Xu} If $k$ is even, then $H'(n,k)\cong H''(n,k)$.}
	
	For an arbitrary fixed element $X$ in $\Omega'_n$, the mapping $\rho_X: \Omega'_n \rightarrow \Omega''_n$ defined by $\rho_X(Y)=X \triangle Y$ is an isomorphism between $H'(n,k)$ and $H''(n,k)$.
	
	{\theorem \label{arc_transive} \cite{Hu_Wang_Xu} $H(n,k)$ is arc-transitive.}
	
	{\lemma \label{Hn} \cite{Hu_Wang_Xu} $H_n$ is a subgroup of $Aut(H(n,k))$.}
	
	{\theorem \label{HnSn}\cite{Hu_Wang_Xu} $Aut(H(n,k))$ contains a subgroup $H_nS_n$ with order $2^nn!$.}
	
	{\theorem \label{two_copy}\cite{Hu_Wang_Xu} If $k$ is odd then $H(n,k)$ is bipartite and connected, and if $k$ is even then $H(n,k)$ consists of two isomorphic connected components.}
	
	Here the bipartite set are $\Omega'_n$ and $\Omega''_n$ hence if any automorphism $\varphi$ of $H(n,k)$ fixes the empty set, then $|\varphi(X)|$ has the same parity as $|X|$. 
	
	\section{When $n=2k-1$}
	In this section will calculate the full automorphism group of $H(n,k)$ where $n=2k-1$ and $k\geq 2$.

	Clearly, $\Omega_n$ is an abelian group under the operation $\triangle$. 
	Now $\Omega_n$ is isomorphic to $\mathbb Z_2 ^n$ by the map $\chi: \Omega_n \rightarrow \mathbb Z_2 ^n$ such that $\chi(X)$ is a $n$ tuple with $i$th position $1$ and $j$th position $0$ if $i \in X$ and $j \notin X$, i.e., $\chi(X)$ is the \textit{characteristic vector} of $X$. So we can consider $\Omega_n$ as a vector space over $\mathbb Z_2$, let us denote this vector space by $V$. Given an element $v=(v_1,v_2, \ldots,v_n )$ of $V$, the \textit{Hamming Weight} of $v$ is the number of non-zero entries in $(v_1,v_2, \ldots ,v_n)$ and it is denoted by $wt(v)$. The addition between two vectors $u,v$ in $V$ is the XOR operation $u \oplus v$. Note that XOR corresponds to the symmetric difference of sets.
 Therefore for $X,Y \in \Omega_n$ we have $\vert X \vert=k$ if and only if $wt(\chi(X))=k$ and $\vert X \triangle Y\vert=k$ if and only if $wt(\chi(X) \oplus \chi(Y))=k$. Let $\sigma \in S_{n+1}$ such that $t=\sigma^{-1}(n+1)$. Define a linear map $T_{\sigma} : V \rightarrow V$ such that 	
	$$T_{\sigma}(v_1,v_2, \ldots,v_n)=
	(v_{\sigma^{-1}(1)}+v_t,v_{\sigma^{-1}(2)}+v_t, \ldots,v_{\sigma^{-1}(n)}+v_t),\mbox{~where~}v_{n+1}=0.$$

	{\lemma \label{bijection}$T_{\sigma}$ is a bijective linear transformation.}
	
	\pf Let $e_i=(0,\ldots,1, \ldots,0)$, where $1$ is in the $i$th position, $\beta_1=\{e_1, e_2, \ldots ,e_n \}$ and $\beta_2=\{e_1, \ldots,e_{t-1},e_{t+1}, \ldots ,e_n,(1,1,\ldots,1)\}$. Clearly $\beta_1$ and  $\beta_2$ are bases of $V$. Therefore we have either $T_{\sigma}(\beta_1)=\beta_1$(set wise, may not maintain the order)  when $t=n+1$ or  $T_{\sigma}(\beta_1)=\beta_2$ (set wise, may not maintain the order) when $t \neq n+1$, hence $T_{\sigma}$ is a bijective linear transformation. \qed

	{\lemma \label{k_preserve} If $wt(v)=k$ then $wt(T_{\sigma}(v))=k$.}
	
	\pf: Let $\sigma \in S_{n+1}$, $t=\sigma^{-1}(n+1)$. If $\sigma(n+1)=n+1$, i.e., $t=n+1$ then $v_t=0$,  hence the result follows. If  $v_t=1$, then $wt(T_{\sigma}(v))=n-k+1=k$, this completes the proof. \qed
	
	Therefore from Lemmas \ref{bijection} and \ref{k_preserve}  we have if $n=2k-1$ then $\vert X \triangle Y \vert=k$ if and only if $wt(\chi(X) \oplus \chi(Y))=k$, hence $wt(T_{\sigma}(\chi(X) \oplus \chi(Y)))=k$, i.e.,  $wt((T_{\sigma} \circ \chi)(X)\oplus (T_{\sigma} \circ \chi)(Y))=k$. So $\chi^{-1}\circ T_{\sigma} \circ\chi$ is an automorphism of $H(2k-1,k)$, call it $f_{\sigma}$.  Therefore   $f_{\sigma}: \Omega_n \rightarrow \Omega_n$ such that
	\begin{equation}\label{f_sigma}
	f_{\sigma }(X)=\left\lbrace \begin{array}{ll}
	\sigma(X) &  \mbox{~when~}  t \notin X\\
	
	(\sigma(X\setminus  \{t\} ))^c & \mbox{~when~}  t \in X.
	\end{array} \right. 
	\end{equation}
	where $t=\sigma^{-1}(n+1)$.
	
	Let $\Gamma=Aut(H(2k-1,k))$ and $\Gamma_{n+1}=\{ f_{\sigma} : \sigma \in S_{n+1} \}$. From \textbf{Appendix A} we have $\Gamma_{n+1}$ is a subgroup of $\Gamma$.  Let $\emptyset$ be the null set and $\Gamma_{\emptyset}$ be the stabilizer of the element $\emptyset$ in $\Gamma$, hence $\Gamma_{n+1} $ is a subgroup of $\Gamma_{\emptyset}$.   As $\emptyset$ is adjacent to only with the sets of size $k$, hence for all $f \in \Gamma_{\emptyset}$ we have $f(\Omega_n ^k)=\Omega_n ^k$.

	{\theorem \label{isom to S_n+1} $\Gamma_{n+1}$ is isomorphic to $S_{n+1}$.}
	\begin{proof}
		Define $\psi: S_{n+1} \rightarrow \Gamma_{n+1}$ such that $\psi(\sigma)=f_{\sigma}$. From \textbf{Appendix A} we have $f_{\sigma_1 \circ \sigma_2}=f_{\sigma_1} \circ f_{\sigma_2}$, hence $\psi(\sigma_1 \circ \sigma_2)=\psi(\sigma_1) \circ \psi(\sigma_2)$, i.e., $\psi$ is a group homomrphism. Clearly $\psi$ is surjective. Let $\sigma \in Ker\psi$ and $\sigma(t)=n+1$. If $t \neq A$ then $f_{\sigma}(\{t\})=A$, which is a contradiction. If $t=n+1$ then $f_{\sigma}(X)=\sigma(X)=X$ for all $X \in \Omega_n$, hence $\sigma=(1)$. So $Ker\psi=\{(1)\}$, i.e., $\psi$ is injective. This completes the proof.                             
	\end{proof}
	
	Note that from the adjacency condition we can conclude that,  irrespective of any parity of $k$,  all neighbours of any set of size $k$ are always even order.  Conversely, every even order set must have neighbours in $\Omega_n^k$.
	
	{\lemma \label{even_image_size} Let $f \in \Gamma_{\emptyset}$. If $\vert X \vert =2p$, $1 \leq p \leq k-1$ then either $\vert f(X) \vert=2p$ or $n-2p+1$.}
	
	\begin{proof}
		The number of neighbours of $X$ in $\Omega_n ^k$ is $\binom{2p}{p} \binom{n-2p}{k-p}= \binom{2p}{p} \binom{2k-2p-1}{k-p}$. Let $u_m=\binom{2m}{m} \binom{2k-2m-1}{k-m}$ be a finite sequence where $1 \leq m \leq k-1$. It is easy to show that $u_m=u_{k-m}$, i.e., number of neighbours of the sets
		of size $2m$ and $2(k-m)=n-2m+1$  are same in $\Omega_n ^k$. Let $1 \leq m<m+1 \leq \frac{k-1}{2}$, i.e., $k-2m-1 \geq 2$. So we have,  $$\frac{u_m}{u_{m+1}}=\frac{(2km+k-2m^2-m)+(k-2m-1)}{2km+k-2m^2-m}.$$ Hence  $\{u_m\}$ is strictly decreasing for $1 \leq m < \frac{k-1}{2}$ 
		and $\{u_m\}$ is strictly increasing for $\frac{k-1}{2} < m \leq k-1$. If $k$ is odd and  $m=\frac{k-1}{2}$ then $u_m=u_{m+1}$, i.e., $u_{\frac{k-1}{2}}=u_{\frac{k+1}{2}}$.  So $k-p$ is the only point that takes the same value as $p$ in the sequence $\{u_m\}$. As graph isomorphism preserves the number of neighbours, hence we have $|f(X) |=2p$ or $n-2p+1$. This completes the proof.
	\end{proof}

	Note that, in previous Lemma if we take $2p=k-1$ then $n-2p+1=k+1$ and if $2p=k+1$ then $n-2p+1=k-1$. Hence if $|X|=k-1$ or $k+1$ then $|f(X)|=k-1$ or $k+1$. If $|X|\neq k-1$ and $k+1$ then $|f(X)| \neq k-1$ and $k+1$.  This result will be used in next Lemma. 
	
	{\lemma \label{one_image}If $\vert X \vert=1$ or $n$ then $\vert f(X) \vert =1$ or $n$ for all $f \in \Gamma_{\emptyset}$.}
	
	\begin{proof}
		If $\vert X \vert=1$, then $X$ has neighbours of size either $k-1$ or $k+1$, and if $\vert X \vert=n$ then $X$ has only neighbours of size  $k-1$. Clearly $|f(X)|$ is odd and $\neq k$. 
		If $1< \vert f(X) \vert  < k$ or $k< \vert f(X) \vert <n$ then $f(X)$ must have  neighbours of size $k-\vert f(X) \vert$ or $ \vert f(X) \vert -k$ respectively and both $k-\vert f(X) \vert$ and $ \vert f(X) \vert -k$ are strictly less than $k-1$. As every neighbour of $X$ will be mapped to every neighbour of $f(X)$, hence $\vert f(X) \vert =1$ or $n$. 
	\end{proof}
	
	{\note \label{note0} If $k$ is even, then from the Theorem \ref{even k} we have $H'(n,k) \cong H''(n,k)$. Consider the isomorphism $\rho_{[n]}$ from $\Omega'_n$ to $\Omega''_n$, hence $\rho_{[n]}(X)=[n] \triangle X= X^c$.}

	{\theorem \label{Size_preserve} Let $f \in \Gamma_{\emptyset}$. If $ f(\{ i \}) =\{ x_i \}$ for all $i \in [n]$, then $f(X)  = \sigma (X)$ for all $X \in \Omega_n$, where $\sigma \in S_{n+1}$ such that $\sigma(i)=x_i \in [n]$ for all $i \in [n]$ and $\sigma(n+1)=n+1$, i.e., $f =f_{\sigma}$ as in the Equation \ref{f_sigma}.}
	
	\begin{proof}
		Let $X\in \Omega_n$ such that $|X|=k-1$. If $i \notin X$ then $\{i \} \sim X$, i.e., $\{x_i\} \sim f(X)$, i.e., $x_i \notin f(X)$ for all $i \notin X$. As $|f(X)|=k-1$ or $k+1$ hence $f(X)=\sigma(X)$ for all $X \in \Omega^{k-1}_n$.
		
		Let $X \in \Omega'_n$ such that $|X|=2p+1$ for $1\leq p \leq k-1$, then  $X$ must have  neighbours of size $k-1$. Consider the set $X'$ which contains $p$ many elements from $X$ and $k-p-1$ many elements from $X^c$. Therefore $|X'|=k-1$ and $X \sim X'$ and  $f(X')=\sigma(X')$. As $f(X)$ is common neighbour of all such $f(X')$, hence $f(X)=\sigma(X)$. Therefore $f(X)=\sigma(X)$ for all $X \in \Omega'_n$. Note that $f([n])=[n]$.
		
		\textbf{Case 1:} Let $k$ be even integer. We proved that $f \in \Gamma_{\emptyset}$ implies $f(X)=\sigma(X)$ for all $X \in \Omega'_n$. From the Theorem \ref{even k} we have $H'(n,k) \cong H''(n,k)$. Therefore from the note \ref{note0} we have  $f(X)=\sigma(X)$ for all $X \in \Omega''_n$, i.e., $f(X)=\sigma(X)$ for all $X \in \Omega_n$.

		\textbf{Case 2:} Let $k$ be odd integer and $X \in \Omega_n^k$. Therefore $X$ must have neighbours of size $k-1$. Let consider the set $X'$ which contains $\frac{k-1}{2}$ many elements from $X$ and $\frac{k-1}{2}$ many elements from $X^c$. Hence $|X'|=k-1$ and $X \sim X'$ and $f(X')=\sigma(X')$. As $f(X)$ is  common neighbour of all such $f(X')$, hence $f(X)=\sigma(X)$ for all $X \in \Omega_n^k$.

		Now let $X \in \Omega''_n$ such that $|X|=2p$ for $1\leq p \leq k$. Therefore $X$ must have neighbours of size $k$.  Let consider the set $X'$ which contains $p$ many elements from $X$ and $k-p$ many from $X^c$. Hence $|X'|=k$ and $X \sim X'$ and $f(X')=\sigma(X')$. As $f(X)$ is  common neighbour of all such $f(X')$, hence $f(X)=\sigma(X)$ for all  $X \in \Omega''_n$. 
		Hence $f(X)=\sigma(X)$ for all $X \in \Omega_n$.
	\end{proof}

	{\note \label{note_1}  Let $\varphi \in \Gamma_{\emptyset}$.  If $\varphi(\{ i \}) =\{ x_i \}$ then by Theorem \ref{Size_preserve} we have $\varphi=f_{\sigma} \in \Gamma_{n+1}$ where $\sigma(n+1)=n+1$. Now consider  $\varphi(\{t\})=[n]$ for some $t \in [n]$. Without lose of generality we can take $t=1$, i.e., $\varphi(\{1\})=[n]$. Then by the Lemma \ref{one_image} we have $ \vert \varphi([n]) \vert=1 $ and $\vert  \varphi(\{j\})\vert =1$ for all $j \in \{2,3, \ldots ,n\}$. 
		
		Let $\tau \in S_{n+1}$ such that $\tau(1)=n+1$, therefore there exists $\sigma \in S_{n+1}$ with $\sigma(n+1)=n+1$ such that $\tau=\sigma \circ (1~n+1)$. Similarly for $\varphi \in \Gamma_{\emptyset}$ with $\varphi(\{1\})=[n]$  there exist
		$f \in \Gamma_{\emptyset}$ such that  $f(\{1\})=[n]$, $f([n])=\{1\}$ and $f(\{i\})=\{i\}$ for all $i \neq 1$ and $\sigma \in S_{n+1}$ with $\sigma(n+1)=n+1$. So we can consider $\varphi([n])=(f_{\sigma} \circ f)([n])$, $\varphi(\{i\})=(f_{\sigma} \circ f)(\{i\})$, where $f_{\sigma}(X)=\sigma(X)$. Therefore $(\varphi \circ f^{-1})(\{i\})=f_{\sigma}(\{i\})=\sigma(\{i\})$, hence by the Theorem \ref{Size_preserve} we have $(\varphi \circ f^{-1})(X)=f_{\sigma}(X)=\sigma(X)$, i.e., $\varphi(X)=(f_{\sigma} \circ f)(X)$ for all $X \in \Omega_n$.
		We will prove that $f=f_{\sigma'}$ as defined in the Equation \ref{f_sigma}, where $\sigma'=(1~ n+1) \in S_{n+1}$ and hence we will have $\varphi=f_{\sigma} \circ f_{\sigma'} \in \Gamma_{n+1}$.
		Now on-wards we will focus on the automorphism $f$. As $f(\{i\})=\{i\}$ for all $i \neq 1$ so by the Theorem \ref{Size_preserve} we have $f(X)=X$ for all $X \in \Omega_n$ not containing $1$.}


	{\lemma \label{k-1_image}  Let $X\in \Omega_n $ such that $|X|=k-1$. If $1 \in X$ then $ f(X) = (X \setminus \lbrace 1 \rbrace)^c $.}
	
	\pf As $|X|=k-1$, so  $X \sim \{i\}$ for all $i \notin X$ and $X \sim [n]$. Hence $f(X) \sim f(\{i\})=\{i\}$ for all $i \notin X$ and $f(X) \sim f([n])=\{1\}$. Therefore $f(X)= (X \setminus \lbrace 1 \rbrace)^c$ is the only option. \qed

	{\lemma \label{odd_image_2}  Let $X\in \Omega'_n $. If $1 \in X$ then $ f(X)  = (X \setminus \lbrace 1 \rbrace)^c$.}
	
	\pf Let $\vert X \vert=2p+1$ where $1 \leq p \leq k-1$. So $X$ must have neighbours of size $k-1$. Consider the set $X'$ which contains $p$ many elements from $X$ including $1$ and $k-p-1$ many elements from $X^c$. Therefore $|X'|=k-1$ and $X \sim X'$, hence $f(X) \sim f(X')$. By previous Lemma \ref{k-1_image} we have $f(X')=(X' \setminus \lbrace 1 \rbrace)^c$. 
	
	Now  Consider the set $X''$ which contains $p$ many elements from $X$ excluding $1$ and $k-p-1$ many elements from $X^c$. Therefore $|X''|=k-1$ and $X \sim X'$, hence $f(X) \sim f(X'')=X''$. 
	
	As $f(X)$ is common neighbour of all such $f(X')$ and $f(X'')$, hence $f(X)= (X \setminus \lbrace 1 \rbrace)^c$. This completes the proof.\qed

	{\lemma \label{even_image_3} Let $X\in \Omega''_n $. If $1 \in X$ then $ f(X)  = (X \setminus \lbrace 1 \rbrace)^c$.}
	
	\pf \textbf{Case 1:} If $k$ is even then by Note \ref{note0} we have $H'(n,k)\cong H''(n,k)$. Hence by previous Lemma \ref{odd_image_2} we have  $f(X)= (X \setminus \lbrace 1 \rbrace)^c$ for all $X \in \Omega''_n$. 
	
	\textbf{Case 2:} Let $k$ be an odd integer. Let $|X|=2p$ where $1 \leq p \leq k-1$. So $X$ must have neighbours of size $k$. Consider the set $X'$ which contains $p$ many elements from $X$ including $1$ and $k-p$ many elements from $X^c$. Therefore $|X'|=k$ and $X \sim X'$, hence $f(X) \sim f(X')$. By previous Lemma \ref{odd_image_2} we have $f(X')=(X' \setminus \lbrace 1 \rbrace)^c$. 
	
	Now  Consider the set $X''$ which contains $p$ many elements from $X$ excluding $1$ and $k-p$ many elements from $X^c$. Therefore $|X''|=k-1$ and $X \sim X'$, hence $f(X) \sim f(X'')=X''$. 
	
	As $f(X)$ is common neighbour of all such $f(X')$ and $f(X'')$, hence $f(X)= (X \setminus \lbrace 1 \rbrace)^c$. This completes the proof.\qed

	{	\note  \label{note_2}  Combining Lemma \ref{even_image_3} and \ref{odd_image_2} we have $f=f_{\sigma'}$ as defined in the Equation \ref{f_sigma},  where $\sigma' \in S_{n+1}$ such that $\sigma'=(1~n+1)$.}
	
	{\theorem \label{subgroup} $\Gamma$ contains a subgroup $H_n\Gamma_{n+1}$   with order $2^n(n+1)!$.}
	
	\pf We already have that the groups $H_n$ and $\Gamma_{n+1}$ are subgroups of $\Gamma$. Let $f_{\sigma} \in \Gamma_{n+1}$ and $\rho_X \in H_n$. Clearly $f_{\sigma}^{-1}=f_{\sigma^{-1}}$. Let $S \in \Omega_n$. $f_{\sigma} \rho_X f_{\sigma}^{-1}(S)=f_{\sigma}\rho_X(f_{\sigma^{-1}}(S))=f_{\sigma}(X \triangle f_{\sigma^{-1}}(S))=\sigma(X) \triangle S=\rho_{\sigma(X)}(S)$. Hence $f_{\sigma} \rho_X f_{\sigma}^{-1}=\rho_{\sigma(X)} \in H_n$. Therefore $H_n\Gamma_{n+1}$ is a subgroup of $\Gamma$. Since $H_n \cap \Gamma_{n+1}$ is the trivial subgroup and   by the Theorem \ref{isom to S_n+1} we have $|\Gamma_{n+1}|=(n+1)!$, hence $|H_n\Gamma_{n+1}|=2^n(n+1)!$. \qed

	{\theorem $Aut(H(2k-1,k))=\Gamma =$ $H_n\Gamma_{n+1}$.}
	
	\pf We have already proved that $H_n\Gamma_{n+1}$ is a subgroup of $\Gamma$ of order $2^n(n+1)!$, hence $\vert \Gamma \vert \geq 2^n(n+1)!$.  We will use the orbit stabilizer Theorem to prove the equality. By the Theorem \ref{arc_transive} we have $H(n,k)$ is vertex-transitive, i.e., the orbit of the null set $\emptyset$ is $\Omega_n$.  Therefore $\vert \Gamma \vert =2^n \vert \Gamma_{\emptyset} \vert $. Clearly $\Gamma_{n+1} $ is a subgroup of $\Gamma_{\emptyset}$.  Let $\varphi \in \Gamma_{\emptyset}$. Then $\varphi(\Omega_n ^k)=\Omega_n^k$. Now by the Lemma  \ref{one_image} we have  $\vert \varphi(\{i\}) \vert =1$ or $n$ for all $i \in [n]$. If $\vert \varphi(\{i\}) \vert =1$ for all $i \in [n]$ then by Theorem \ref{Size_preserve} we have $\varphi(X)  = \sigma (X)$ for all $X \in \Omega_n$, where $\sigma \in S_{n+1}$ and $\sigma(n+1)=(n+1)$, i.e., $\varphi \in \Gamma_{n+1}$. 
	
	Now let  $\varphi(\{i\})=[n]$ for some $i \in [n]$. Then by Note \ref{note_1} and \ref{note_2} we have $\varphi=\sigma \circ f_{\sigma'}$ for some $\sigma, \sigma' \in S_{n+1}$ with $\sigma(n+1)=n+1$ and $\sigma'=(1~n+1)$. Hence $\varphi \in \Gamma_{n+1}$. As $\varphi$ is arbitrary element from $\Gamma_{\emptyset}$, hence $\Gamma_{\emptyset} \subseteq \Gamma_{n+1}$ and hence $\Gamma_{\emptyset}= \Gamma_{n+1}$. Therefore $\vert \Gamma \vert =2^n \vert \Gamma_{n+1} \vert =2^n(n+1)!$, hence $Aut(H(2k-1,k))=\Gamma =H_n\Gamma_{n+1}$.  This completes the proof. \qed

	\section{When $k$ is an odd integer and $n=2k+1$}
	
	In this section, we will investigate the full automorphism group of $H(n,k)$ where
	$k$ is an odd integer and  $n=2k+1$. Let $\Gamma'=Aut(H(2k+1,k))$ and $\Gamma'_{\emptyset}$ be the stabilizer of null set in $\Gamma'$. Therefore if $f \in \Gamma'_{\emptyset}$ then $f(\Omega_n^k)=(\Omega_n^k).$
	{\lemma \label{even_image_size2} Let $f \in \Gamma'_{\emptyset}$. If $\vert X \vert =2p$, $1 \leq p \leq k$ then either $\vert f(X) \vert=2p$ or $n-2p+1$.}
	
	\begin{proof}
		
		Note that as $|X|\leq 2k$, hence $X$ must have neighbours in $\Omega^k_n$. The number of neighbours of $X$ in $\Omega_n ^k$ is $\binom{2p}{p} \binom{n-2p}{k-p}= \binom{2p}{p} \binom{2k-2p+1}{k-p}$. Let $u_m=\binom{2m}{m} \binom{2k-2m+1}{k-m}$ be a finite sequence where $1 \leq m \leq k$. It is easy to show that $u_m=u_{k-m+1}$, i.e., number of neighbours of the sets
		of size $2m$ and $2(k-m+1)=n-2m+1$  are same in $\Omega_n ^k$. 
		$$\frac{u_m}{u_{m+1}}=\frac{(2km+k-2m^2+m+1)+(k-2m)}{2km+k-2m^2+m+1}.$$ Hence  $\{u_m\}$ is strictly decreasing for $1 \leq m < \frac{k}{2}$ 
		and $\{u_m\}$ is strictly increasing for $\frac{k}{2} < m \leq k$. So $k-p+1$ is the only point that takes the same value as $p$ in the sequence $\{u_m\}$. As graph isomorphism preserves the number of neighbours, hence we have $|f(X) |=2p$ or $n-2p+1$. This completes the proof.
	\end{proof}

	Hence if $|X|=k+1$ then $|f(X)|=n-(k+1)+1=k+1=|X|$ for all $f \in \Gamma'_{\emptyset}$, i.e., $f(\Omega^{k+1}_n)=\Omega^{k+1}_n$. Therefore $f([n])=[n]$ for all $f \in \Gamma'_{\emptyset}$.

	{\lemma \label{odd_image_size2} Let $f \in \Gamma'_{\emptyset}$. If $\vert X \vert =2p-1$, where  $1\leq p \leq k$ then either $\vert f(X) \vert=2p-1$ or $n-2p$.}
	
	\pf Note that as $|X|\leq 2k+1$, hence $X$ must have neighbours of size $k+1$. The number of neighbours of $X$ of size $k+1$ is $\binom{2p-1}{p} \binom{n-(2p-1)}{k-p+1}=\binom{2p-1}{p} \binom{2k+2p+2}{k-p+1}$. Let $u_m=\binom{2m-1}{m} \binom{2k+2m+2}{k-m+1}$ be a finite sequence where $1 \leq m \leq k+1$. It is easy to show that $u_{m}=u_{k-m+1}$, i.e.,  number of neighbours of the sets of size $2m-1$ and $2(k-m+1)-1=n-2m$ in $\Omega_n ^{k+1}$ are same. 
	By some easy computation we can show $$\frac{u_m}{u_{m+1}}=\frac{(2km+k-2m^2+m+1)+(k-2m)}{2km+k-2m^2+m+1}.$$ Hence $\{u_m\}$ is strictly decreasing sequence for $1 \leq m < \frac{k}{2}$ and strictly increasing sequence for $\frac{k}{2} < m \leq k$.  So $k-p+1$ is the only point that takes the same value as $p$ in the sequence $\{u_m\}$. As graph isomorphism preserves the number of neighbours, hence we have $|f(X) |=2p-1$ or $n-2p$. This completes the proof. \qed

	Hence if $|X|=1$ then either $|f(X)|=1$ or $n-2$ for all $f \in \Gamma'_{\emptyset}$.
	
	{\theorem \label{Size_preserve2} Let $f \in \Gamma'_{\emptyset}$. If $ f(\{ i \}) =\{ x_i \}$ for all $i \in [n]$, then $f(X)  = \sigma (X)$ for all $X \in \Omega_n$, where $\sigma \in S_{n}$ such that $\sigma(i)=x_i$ for all $i \in [n]$.}
	
	The proof is similar to that of Theorem \ref{Size_preserve}.
	
	{\lemma \label{singleton_image} Let $f \in \Gamma'_{\emptyset} $. If $|f(\{i\})|=n-2$ for some $i \in [n]$ then the following conditions hold.
		\begin{itemize}
			\item There exists $t \in[n]$ such that $|f(\{t\})|=1 $.
			\item $|f(\{i\})|=n-2$ for all $i \in [n]\setminus \{t\}$.
			\item  $f(\{i\})\cap f(\{t\})=\emptyset$ for all $i \in[n]\setminus \{t\}$.
			\item $|f(\{i\}) \cap f(\{j\})|=n-3$ for all $i\neq j\in [n]\setminus \{t\}$. 
		\end{itemize}
	}
	
	\pf Let us consider $i \neq j$ with  $|f(\{i\})|=n-2=|f(\{j\})|$. $|f(\{i\}) \cap f(\{j\})|=|f(\{i\})|+|f(\{j\})|-|f(\{i\}) \cup f(\{j\})|\geq (n-2)+(n-2)-n=n-4$. Therefore $|f(\{i\}) \cap f(\{j\})|=n-4$ or $n-3$. Let $|f(\{i\}) \cap f(\{j\})|=n-4$. 
	The number of common neighbours of $\{i\}$ and $\{j\}$ in $\Omega_n^{k+1}$ is $\binom{2k-1}{k-1}$ which is not equal to the number of common neighbours of  $f(\{i\})$ and $f( \{j\})$ in $\Omega_n^{k+1}$ which is $\binom{2k-3}{k-1}\binom{2}{1}\binom{2}{1}$. 
	Hence $|f(\{i\}) \cap f(\{j\})|=n-3$. 
	Note that there are at most $(n-1)$ many such sets in $\Omega_n$.
	Hence if $|f(\{i\})|=n-2$ for some $i \in [n]$ then there exists at least one $t$ such that $|f(\{t\})|=1$. 
	
	If possible let $f(\{t\}) \in f(\{i\})$. The number of common neighbours of $\{i\}$ and $\{t\}$ in $\Omega_n^{k+1}$ is $\binom{2k-1}{k-1}$,  which is not equal to the number of common neighbours of $f(\{i\})$ and $f( \{t\})$ in $\Omega_n^{k+1}$ which is $\binom{2k-2}{k-1}\binom{2}{1}$, which is a contradiction. Hence $f(\{t\}) \notin f(\{i\})$. Therefore there exists exactly one $t$ such that $|f(\{t\})|=1$ and $|f(\{i\})|=n-2$ for all $i \in [n]\setminus \{t\}$.  This completes the proof. \qed

	Form the above lemma we have if $|f(\{i\})|=n-2$ for some $i \in [n]$ then $|f(\{t\})|=1$ and $|f(\{i\})|=n-2$ for all $i \in [n]\setminus \{t\}$.  $|f(\{i_1\}) \cap f(\{i_2\})|=n-3=(n-2)-1$, $|f(\{i_1\}) \cap f(\{i_2\}) \cap f(\{i_3\})|=(n-2)-2$, \ldots , $|f(\{i_1\}) \cap f(\{i_2\}) \cap \cdots \cap f(\{i_r\})|=(n-2)-(r-1)=n-r-1$ for all $i_j\in [n]\setminus \{t\}$, $j=1,2,\ldots ,r$.
	
	Let $\sigma \in S_{n+1}$ such that $\sigma(t)=n+1$ where $t \in [n]$. Define $f^{\sigma}: \Omega_n \rightarrow \Omega_n$ such that $f^{\sigma}(\{t\})=\{\sigma(n+1)\}= \{ t'\}$ and of course $ f^{\sigma}([n])=[n]$. Now we have $n-1$ many subsets of $S=[n]\setminus \{t'\}$ of size $n-2$ such that size of the intersection of any two subsets is $n-3$.  
	Take $f^{\sigma}(\{i\})=S\setminus \{\sigma(i)\}$ for all $i \in [n]\setminus \{t\}$ and


	
	
	

	\begin{equation}\label{f_sigma4}
	f^{\sigma }(X)=\left\lbrace \begin{array}{ll}
	\sigma(X) &  \mbox{~if~} |X| \mbox{~is even  and~}t \notin X\\
	(\sigma(X \setminus \{t\}))^c & \mbox{~if~} |X| \mbox{~is even  and~}t \in X\\
	(\sigma(X \cup \{n+1\}))^c & \mbox{~if~} |X| \mbox{~is odd  and~}t \notin X\\
	\sigma(X\cup \{n+1\} \setminus \{t\}) &  \mbox{~if~} |X| \mbox{~is odd  and~}t \in X.
	\end{array} \right. 
	\end{equation} 
	
	If $\sigma(n+1)=n+1$ then we take $f^{\sigma}(X)=\sigma(X)$ for all $X \in \Omega_n$.
	Let $\Gamma^{n+1}=\{ f^{\sigma} : \sigma \in S_{n+1} \}$. As $k$ is odd integer, hence by  some basic set theoretical calculations we can show that $|X \triangle Y|=k \Leftrightarrow |f^{\sigma}(X) \triangle f^{\sigma}(Y)|=k$. Clearly $\Gamma^{n+1} \subseteq \Gamma'_{\emptyset}$.
	
	By some coputations we can show that $\Gamma^{n+1}$ is a subgroup of $\Gamma'$ ( see \textbf{Appendix B}). 
	

	{\theorem  $\Gamma^{n+1}$ is isomorphic to $S_{n+1}$.}
	
	The proof is similar as that of Theorem \ref{isom to S_n+1}.

	{\theorem \label{f_sigma=varphi} Let $f \in \Gamma'_{\emptyset}$. If $|f(\{i\})|=n-2$ for some $i \in \{1,2,\ldots ,n\}$ then $f =f_{\sigma}$ as defined in the Equation \ref{f_sigma4}.} 
	
	\pf As $|f(\{i\})|=n-2$ for some $i \in \{1,2,\ldots ,n\}$ so $f$ must satisfy the Lemma \ref{singleton_image}. Let $t \in  \{1,2,\ldots ,n\}$ such that $|f(\{t\})|=1$. 
	
	\textbf{Case 1:} Let $|X|=k+1$ and $t \notin X$. Let $X=\{x_1, \ldots ,x_{k+1}\}$. So $\{x_1\}, \ldots , \{x_{k+1}\} \sim X$ imply $f(\{x_i\}) \sim f(X)$ for all $i=1,2, \ldots , k+1$, i.e., $|S \setminus \{ \sigma(x_i) \} \cap f(X)|=k$, i.e., $|(S \setminus \{ \sigma(x_i) \} )^c \cap f(X)|=1$, i.e., $|(S^c \cup \sigma(x_i)) \cap f(X)|=1$ for all  $i=1,2, \ldots , k+1$. Clearly $t'(=\sigma(n+1)) \notin f(X)$, hence $\sigma(x_i) \in f(X)$ for all $i=1,2, \ldots , k+1$, i.e., $f(X)=\sigma(X)$.

	\textbf{Case 2:} Let $|X|=k+1$ and $t \in X$.  Let $X=\{t,x_1, \ldots ,x_{k}\}$. So $\{t\}, \{x_1\}, \ldots , \{x_{k}\} \sim X$ imply $f(\{t \}) \sim f(X)$ and   $f(\{x_i\}) \sim f(X)$ for all $i=1,2, \ldots , k$ , i.e.,  $|S \setminus \{ \sigma(x_i) \} \cap f(X)|=k$ for all $i=1,2, \ldots , k$. If possible let $\sigma(x_i) \in f(X)$ for some $i \in \{1,2,\ldots , k\}.$ Therefore  $|S \setminus \{ \sigma(x_i) \} \cap f(X)|=k$ implies  $|S  \cap f(X)|=k+1$, i.e., $f(X) \subset S$. But $f(\{t \})=\{t'\} \sim f(X)$ implies $t ' \in f(X)$, which is a contradiction as $t' \notin S$. Hence $f(X)=(\sigma(X \setminus \{t\}))^c$.

	\textbf{Case 3:} Let $|X|=2p-1$, $1 \leq p \leq k$ and $t \notin X$. Let $X=\{ x_1, \ldots ,x_{2p-1}\}$ and $X^c=\{t, y_1, y_2, \ldots , y_{n-2p} \}$. Consider the set $X'=\{x_1, \ldots , x_p, t, y_1, \ldots , y_{k-p} \}.$ Hence $|X'|=k+1$ and $X' \sim X$, i.e., $f(X') \sim f(X)$, i.e., $(\sigma(X'\setminus \{t\}))^c \sim f(X)$. As $|f(X)|=2p-1$ or $n-2p$, hence    $|(\sigma(X'\setminus \{t\}))^c \cap f(X)|=p$ or $(k-p+1)$, i.e., $|\{ \sigma(x_{p+1}), \ldots , \sigma(x_{2p-1}), \sigma(y_{k-p+1}), \ldots , \sigma(y_{n-2p}) \} \cap f(X)|=p$ or $(k-p+1)$. We can vary $x_1, \ldots , x_{p}$ from $X$, hence $\sigma(y) \in f(X)$ for all $y \in (X \cup \{n+1\})^c$, i.e., $|f(X)| \geq n-2p$. If $ n-2p \geq 2p$, then the result follows.
	
	Let $n-2p < 2p$.  If possible let $|f(X)|=2p$, i.e.,  $f(X)$ contains remaining $(4p-n)$ many elements from $\sigma(X)$. Let $f(X) \cap \sigma(X)=\{ \sigma(x_1), \ldots , \sigma(x_{4p-n})  \}.$ If $4p-n<p$ then consider the set $X''=\{ x_{2p}, \ldots , x_{4p-n+1}, x_{4p-n}, \ldots , x_{p+1},  y_1, \ldots , y_{k-p} \}$. If $p \geq 4p-n$ then consider the set $X''=\{ x_2, \ldots , x_{p+1}, y_1, \ldots , y_{k-p}  \}$. In both of the cases we have $X \sim X''$, i.e., $f(X) \sim f(X'')$ but $(k-p+1) \lneq |f(X) \cap f(X'')| \lneq p$, which is the contradiction. Hence $|f(X)|=n-2p$, i.e., $f(X)=(\sigma(X \cup \{n+1\}))^c$.

	\textbf{Case 4:} Let $|X|=2p-1$, $1 \leq p \leq k$ and $t \in X$.   Let $X=\{ t, x_2, \ldots , x_{2p-1} \}$. Consider the set $X'=\{ x_2, \ldots , x_{p+1}, y_1, \ldots , y_{k-p+1} \}$, where $y_1, \ldots , y_{k-p+1} \in (X \cup \{t\})^c$.  Hence $|X'|=k+1$ and $X' \sim X$, i.e., $f(X') \sim f(X)$, i.e., $\sigma(X') \sim f(X)$, As $|f(X)|=2p$ or $n-2p$, hence $|\sigma(X') \cap f(X)|=p$ or $(k-p+1)$, i.e., $| \{\sigma(x_2), \ldots , \sigma(x_{p+1}), \sigma(y_1), \ldots , \sigma(y_{k-p+1})\} \cap f(X) |=p$ or $(k-p+1)$. As both $\{ x_2, \ldots ,x_{p+1} \}$ and $\{ y_1, \ldots , y_{k-p+1} \}$ are arbitrary from $X$ and $X^c$ respectively, so we have two options either $\sigma(X \setminus \{t\}) \subset f(X)$, i.e., $|f(X)| > 2p-2$ or $\sigma(X^c) \subset f(X)$, i.e., $|f(X)|>n-2p+1$. As $|f(X)|=2p-1$ or $n-2p$, hence $2p-1<n-2p+1$ implies $\sigma(X \setminus \{t\}) \subset f(X)$.
	
	Now let $n-2p+1<2p-1$. If possible let $\sigma(X^c) \subset f(X)$, hence remaining $(4p-n-2)$ many elements of $f(X)$ are from $\sigma(X \setminus \{t\})$. Let $f(X) \cap \sigma(X \setminus \{t\})=\{ x_2, \ldots , x_{4p-n-1} \}.$  
	If $4p-n-2 <p$ then consider the set $X''=\{ x_{2p}, \ldots , x_{4p-n+1}, x_{4p-n}, \ldots , x_{p+1}, y_1, \ldots , y_{k-p+1} \}$. If $p \geq 4p-n$ then consider the set $X''=\{ x_3, \ldots , x_{p+2}, y_1, \ldots , y_{k-p+1}  \}$. In both of the cases we have $X \sim X''$, i.e., $f(X) \sim f(X'')$ but $(k-p+1) \lneq |f(X) \cap f(X'')| \lneq p$, which is the contradiction. Hence  $\sigma(X^c) \not\subset f(X)$, i.e., $\sigma(X \setminus \{t\}) \subset f(X)$. By choosing suitable set we can  prove that $\sigma(y) \notin f(X)$ for all $y \in X^c$. Note that $t' \notin \sigma(X^c)$. As $|\sigma(X \setminus \{t\}|=2p-2$ and $|f(X)|=2p-1$, hence $t' \in f(X)$, i.e., $f(X)=\sigma(X \cup \{n+1\} \setminus \{t\})$.
	
	\textbf{Case 5:} Let $|X|=2p$, $1 \leq p \leq k$ and $t \notin X$. Let $X=\{x_1, \ldots , x_{2p}\}$. Consider the set $X'=\{x_1, \ldots , x_p,t,y_2, \ldots , y_{k-p}\}$, where $y_2, y_3, \ldots , y_{k-p} \in (X \cup \{t\})^c$.  Hence $|X'|=k$ and $X' \sim X$, i.e., $f(X') \sim f(X)$, i.e., $\sigma(X \cup \{n+1\} \setminus \{t\}) \sim f(X)$. As $|f(X)|=2p$ or $n-2p+1$, hence $|\sigma(X \cup \{n+1\} \setminus \{t\}) \cap f(X)|=p$ or $(k-p+1)$, i.e., $|\{t',\sigma(x_1), \ldots , \sigma(x_p),\sigma(y_2), \ldots , \sigma(y_{k-p})\} \cap f(X)|=p$ or $(k-p+1)$ where $\sigma(n+1)=t'$.  We can vary $y_2, y_3, \ldots , y_{k-p}$ from $ (X \cup \{t\})^c$. Hence $\sigma(x_1), \sigma(x_2), \ldots , \sigma(x_p) \in f(x)$. As $x_1, \ldots , x_p\}$ are arbitrary from $X$, hence $\sigma(x_i) \in f(X)$ for all $i=1,2, \ldots , 2p.$ , i.e., $|f(X)| \geq 2p$. If $2p \geq n-2p+1$, then we have $|f(X)|=2p$. 
	
	Let $2p<n-2p+1$, i.e., $k< 2p$, i.e., $k-p+1> n-4p+1$. If possible let $|f(X)|=n-2p+1$, i.e.,  $f(X)$ contains remaining $(n-4p+1)$ many elements from $\sigma(X)^c$. Let $f(X) \cap \sigma(X)^c=\{\sigma(y_1),\ldots ,\sigma( y_{n-4p+1})\}.$  Consider the set $X''=\{x_1,x_2, \ldots , x_p, y_2, \ldots ,y_{n-4p+1}, \ldots , y_{k-p}, y_{k-p+1}\}$ where $y_{i} \in  (X \cup \{t\})^c$ for all $i$.   Hence $|X''|=k$ and $X \sim X''$, i.e., $f(X) \sim f(X'')$, i.e., $(\sigma(X \cup \{n+1\}))^c \sim f(X)$, i.e. $|(\sigma(X \cup \{n+1\}))^c \cap  f(X)|=p$ or $(k-p+1)$. But $  (\sigma(X \cup \{n+1\}))^c \cap  f(X)= \{\sigma(x_{p+1}), \ldots , \sigma(x_{2p}), \sigma(y_1) \}$, i.e., $| (\sigma(X \cup \{n+1\}))^c \cap  f(X)|=p+1$,  which is a contradiction.  Hence $|f(X)|=2p$, i.e., $f(X)=\sigma(X)$.

	\textbf{Case 6:} Let $|X|=2p$, $1 \leq p \leq k$  and $t \in X$. Let $X=\{ t, x_2, \ldots , x_{2p} \}$ and $X^c=\{ y_1, y_2, \ldots , y_{n-2p} \}$. Consider the set $X'=\{ x_2, \ldots , x_{p+1}, y_1, \ldots , y_{k-p} \}$.  Hence $|X'|=k$ and $X \sim X'$, i.e., $f(X) \sim f(X')$, i.e., $(\sigma(X \cup \{n+1\}) )^c \sim f(X)$. As $|f(X)|=2p$ or $n-2p+1$,  hence $|(\sigma(X \cup \{n+1\} ))^c \cap f(X)|=p$ or $(k-p+1)$, i.e., $|\{ \sigma(x_{p+2}), \ldots , \sigma(x_{2p}), \sigma(y_{k-p+1}), \ldots , \sigma(y_{n-2p}) \} \cap f(X)|=p$ or $(k-p+1)$. We can vary $x_2, \ldots , x_{p+1}$ from $(X \setminus \{t\})$, hence $\sigma(y) \in f(X)$ for all $y \in X^c$, i.e., $|f(X)| \geq n-2p+1$. If $ n-2p+1 \geq 2p$, then we have $|f(X)|=n-2p+1$. 
	
	Let $n-2p+1 < 2p$.  If possible let $|f(X)|=2p$, i.e.,  $f(X)$ contains remaining $(4p-n-1)$ many elements from $\sigma(X)$. Let $f(X) \cap \sigma(X)=\{ \sigma(x_2), \ldots , \sigma(x_{4p-n})  \}.$ If $4p-n-1 <p$ then consider the set $X''=\{ x_{2p}, \ldots , x_{4p-n+1}, x_{4p-n}, \ldots , x_{p+2}, t, y_1, \ldots , y_{k-p} \}$. If $p \geq 4p-n$ then consider the set $X''=\{ x_3, \ldots , x_{p+1},t, y_1, \ldots , y_{k-p}  \}$. In both of the cases we have $X \sim X''$, i.e., $f(X) \sim f(X'')$ but $(k-p+1) \lneq |f(X) \cap f(X'')| \lneq p$, which is the contradiction. Hence $|f(X)|=n-2p+1$, i.e., $f(X)=(\sigma(X \setminus \{t\}))^c$. 

 Therefore $f =f_{\sigma}$ as defined in the Equation \ref{f_sigma4}. \qed

	{\theorem $Aut(H(2k+1,k))=\Gamma' = H_n\Gamma^{n+1}$.}
	
	\pf Arguing as same as the Theorem \ref{subgroup} we have $\Gamma'$ contains a subgroup $H_n\Gamma^{n+1}$   with order $2^n(n+1)!$.  From the  Note \ref{Size_preserve2} and Theorem \ref{f_sigma=varphi} we have $\Gamma'_{\emptyset}=\Gamma^{n+1}$. Hence by orbit-stabilizer Theorem we have $|Aut(H(2k+1,k))|=|H_n\Gamma^{n+1}|=2^n(n+1)!$, i.e., $Aut(H(2k+1,k))=\Gamma' = H_n\Gamma^{n+1}$. \qed

	\section{When $k$ is an even integer and $n=2k+1$}
	
	Fu-Tao Hu \textit{et.al.} \ref{HnSn} proved that $H_nS_n$ is a subgroup of $Aut(H(n,k))$. In this section, we will prove the equality for the case $n= 2k+1$ when $k$ is an even integer. 

	{\theorem \label{Aut(H(n,k))0} $Aut(H(n,k))=H_nS_n$, where $k$ is an even integer and $n=2k+1$. }
	
	\begin{proof}
		From the Theorem \ref{HnSn} we have  $|Aut(H(n,k))|\geq |H_n S_n|=2^n n!$. As $H(n,k)$ is vertex-transitive, to prove the equality we will use the orbit-stabilizer theorem, i.e., we will show the stabilizer of null set $\emptyset$ equals to the group $S_n$.  Let $f \in Aut(H(n,k))$ such that $f(\emptyset)=\emptyset$. If possible let $f \notin S_n$. Note that if $f \notin S_n$, then $f$ must satisfy the equation \ref{f_sigma4}.
		
		\textbf{Case 1:} Let $\frac{k}{2}$ be odd integer. Consider two sets $X$ and $Y$ from $\Omega'_n$ such that $|X|=|Y|=\frac{k}{2}$, $t \in X$ and  $X \cap Y=\emptyset$. Let $X'=X \setminus \{t\}$.
		Therefore $X \sim Y$ and $f(X)=\sigma(X \cup \{n+1\} \setminus \{t\})=\sigma(X')\cup \{t'\}$, $f(Y)=\sigma(Y \cup \{n+1\})^c$. As $X \cap Y=\emptyset$, hence $\sigma(X') \cap \sigma(Y)=\emptyset$, i.e., $\sigma(X') \subset \sigma(Y)^c\setminus \{t \}$. Therefore $f(X) \triangle f(Y)=(\sigma(Y)^c \setminus \sigma(X')) \cup \{t\}$, i.e., $|f(X) \triangle f(Y)|=k+2$, i.e., $f(X) \nsim f(Y)$, which is a contradiction as $X \sim Y$.
		
		\textbf{Case 2:}  Let $\frac{k}{2}$ be even integer. Consider two sets $X$ and $Y$ from $\Omega''_n$ such that $|X|=|Y|=\frac{k}{2}$, $t \in X$ and  $X \cap Y=\emptyset$. Let $X'=X \setminus \{t\}$. Hence $X \sim Y$. Similarly as previous case we can prove that $f(X) \nsim f(Y)$, which is a contradiction.
		
		Combining two cases we have $f \in S_n$. So the stabilizer of null set $\emptyset$ equals to the group $S_n$. By orbit-stabilizer theorem we have $|Aut(H(n,k))|=2^n n!$ and hence $Aut(H(n,k))=H_nS_n$. 
	\end{proof}

	\section{When $n \neq 2k-1,2k,2k+1$}
	
	We already have that $H_nS_n$ is a subgroup of $Aut(H(n,k))$. In this section, we will prove the equality for the cases $n\neq 2k-1,2k,2k+1$. 
	
	{\theorem \label{Aut(H(n,k))1} $Aut(H(n,k))=H_nS_n$, where $k$ is an odd integer and $n \geq 2k+2$. }
	
	\begin{proof}
		
		From the Theorem \ref{HnSn} we have  $|Aut(H(n,k))|\geq |H_n S_n|=2^n n!$. As $H(n,k)$ is vertex-transitive, to prove the equality we will use the orbit-stabilizer theorem, i.e., we will show the stabilizer of null set $\emptyset$ equals to the group $S_n$.  Let $f \in Aut(H(n,k))$ such that $f(\emptyset)=\emptyset$.

		\begin{figure}[ht]
			\centering
			\begin{center}
				\includegraphics[scale=1.5]{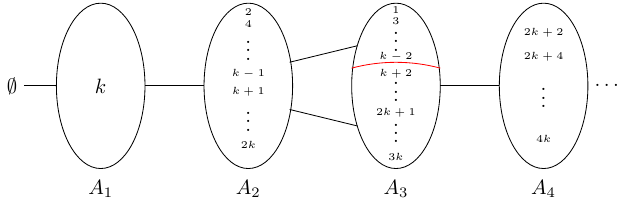}
				\label{pic}
			\end{center}
		\end{figure}

		Let $A_i$ denote the collection of sets which are at distance $i$ from $\emptyset$, hence $A_1=\Omega_n ^k$. So $f(A_i)=A_i$ for all $i$ and $|f(X)|=|X|$ for all $X \in A_1$. Note that as $n \geq 2k+2$, any set from $A_3$ of size $\geq k+2$ must have a neighbour in $A_4$, but the sets of size $\leq k-2$ has no neighbour in $A_4$. Therefore $|f(\{i\})| \leq k-2$. We know that $H(2k-1,k)$ is induced subgraph of $H(n,k)$. So $f \vert_{[2k-1]} \in S$ where $S$ is subgroup of $S_n$ and $S \cong S_{2k-1}$. Hence from the Lemma \ref{one_image} we have $|f(\{i\})|=1$ for all $i=1, \ldots , 2k-1$. 
		
		Let $\{ a\}, \{b\} \in \Omega_n$. Number of common neighbours of size $k-1$ is $\binom{n-2}{k-1}$ and of size $k+1$ is $\binom{n-2}{k-1}$. Therefore total number of common neighbours is $2 \binom{n-2}{k-1}$. Let $X \in A_3$ such that $|X|=2p+1 \leq k-2$. 
		
		\textbf{Case 1:}  Let $a \in X$ and $p \geq 1$. Number of common neighbours of size $k-1$ is $\binom{2p}{p} \binom{n-2p-1}{k-p-1}$ of size $k+1$ is $\binom{2p}{p} \binom{n-2p-1}{k-p}$. Therefore total number of common neighbours is $\binom{2p}{p} \binom{n-2p}{k-p}=u_p$(say). $X$ has neighbours of size $k-(2i+1)$ and $k+(2i+1)$ for all $i \in \{0, \ldots ,p \}$ and total number of neighbours is fixed for every vertex that is $\binom{n}{k}$, so without calculation  we can say $\{u_p\}$ should be a decreasing sequence for $p=1, \ldots , \frac{k-3}{2}$ and also from some easy computations (see \textbf{Appendix C}) we can prove that  $\{ u_p \}$ is a decreasing sequence and most importantly $u_1=2 \binom{n-2}{k-1} \gneq u_i$ for all $i=2, \ldots , \frac{k-3}{2}$. Therefore $|f(\{i\})|=1$ or $3$.

		\textbf{Case 2:}  Let $a \notin X$. Number of common neighbours of size $k-1$ is $\binom{2p+1}{p} \binom{(n-2p-1)-1}{k-p-1}$ of size $k+1$ is $\binom{2p+1}{p+1} \binom{(n-2p-1)-1}{k-p-1}$. Therefore total number of common neighbours is $\binom{2p+2}{p+1} \binom{n-2p-2}{k-p-1}=v_p$(say). Same as case 1  we can prove that  $\{ v_p \}$ is a decreasing sequence and $v_0\gneq v_i$ for all $i=1, \ldots , \frac{k-3}{2}$.
		
		We already have $|f(\{i\})|=1$ for all $i=1, \ldots , 2k-1$. Let $f(\{i\})=\{x_i\}$ for all $i=1, \ldots , 2k-1$. We claim that $|f(\{i\})|=1$ for all $i=1, \ldots , n$.  If not, let $|f(\{a\})|>1$ for some $a \in \{2k, \ldots , n\}.$ As $\{i\}$ and $\{a\}$ have $2 \binom{n-2}{k-1}$ many common neighbours for all $i=1,\ldots , 2k-1.$   so $\{x_i\}$ and $f(\{a\})$ also have  $2 \binom{n-2}{k-1}$ many common neighbours for all $i=1, \ldots , 2k-1$. From above two cases we have either $|f(\{a\})|=3$ or $|f(\{a\})|=1$. If $|f(\{a\})|=3$ then  $x_i \in f(\{a\})$, i.e., $f(\{a\})$ must contain atleast $2k-1$ many elements, which is a contradiction. Therefore  $|f(\{i\})|=1$ for all $i=1, \ldots , n$, hence $f(\{i \})$ has only neighbours of size $k+1$ and $k-1$.

		Let $f(\{i\})=\{ x_i \}$ for all $i \in \{1,2, \ldots , n\}$ and $\sigma \in S_n$ such that $\sigma(i)=x_i$. We will show that $f(X)=\sigma (X)$ for all $X \in \Omega_n$ by induction. Let $X\in \Omega_n^{k+1}$. If $i \in X$ then $\{i \} \sim X$, i.e., $f(\{i \}) \sim f(X)$, i.e., $x_i \in f(X)$ for all $i \in X$. As $|f(X)|=k+1$ so $f(X)=\sigma(X)$ for all $X\in \Omega_n^{k+1}$.
  
  Let $X \in A_1$, i.e., $|X|=k$. Let consider the set $X'$ which contains $\frac{k+1}{2}$ many elements from $X$ and $\frac{k+1}{2}$ many elements from $X^c$. Hence $|X'|=k+1$ and $X \sim X'$ and $f(X')=\sigma(X')$. As $f(X)$ is  common neighbour of all such $f(X')$, hence $f(X)=\sigma(X)$ for all $X \in A_1$, i.e., the statement is true for all $X \in A_1$. 
		
		Let $X \in A_2$, $|X|=2l$.  Let consider the set $X'$ which contains $l$ many elements from $X$ and $k-l$ many from $X^c$. Hence $|X'|=k$ and $X \sim X'$ and $f(X')=\sigma(X')$. As $f(X)$ is  common neighbour of all such $f(X')$, hence $f(X)=\sigma(X)$, i.e,   the statement is true for all $X \in A_2$.

		Let the statement is true for all $X \in A_i$, $1\leq i \leq m$. Let $Y \in A_{m+1}$ where $m+1 \geq 3$ and $|Y|=l$. If $l>k$ then consider $Y_{i_1 i_2\ \ldots i_k}=Y \setminus \{ i_1,i_2,\ldots ,i_k \}$. So $Y\sim Y_{i_1i_2\ldots i_k}$ and $Y_{i_1 i_2\ldots i_k} \in A_m$. Hence by induction hypothesis we have $f(Y_{i_1i_2\ldots i_k})=\sigma(Y_{i_1i_2\ldots i_k})$. As $f(Y)$ is  common neighbour of  $f(Y_{i_1i_2\ldots i_k})$ for all $i_1,\ldots ,i_k \in Y$, hence $f(Y)=\sigma(Y)$. If $l<k$, i.e, $Y \in A_3$ and $1<l \leq k-2$ then consider $Y_{i_1i_2\ldots i_k}=Y \sqcup \{ i_1,i_2,\ldots ,i_k \}$ where $i_1,i_2,\ldots ,i_k \notin Y$. Then $Y_{i_1 i_2\ldots i_k} \in A_2$. Hence similarly we can prove that $f(Y)=\sigma(Y)$.  So by the principle of induction, we have $f(X)=\sigma(X)$ for all $X \in \Omega_n$ and hence $f \in S_n$. So the stabilizer of null set $\emptyset$ equals to the group $S_n$. By orbit-stabilizer theorem we have $|Aut(H(n,k))|=2^n n!$ and hence $Aut(H(n,k))=H_nS_n$. 
	\end{proof}

	{\theorem \label{Aut(H(n,k))2} $Aut(H(n,k))=H_nS_n$, where $k$ is an odd integer and $n \leq 2k-2$. }

	\begin{proof}
		If $n$ is even then from the Theorem \ref{H(n,n-k)} we have $H(n,k) \cong H(n,n-k)$. As $n \leq 2k-2$, so $n-k \geq 2(n-k)+2$. Therefore by previous Theorem we have $Aut(H(n,k))=Aut(H(n,n-k)) =H_nS_n$.
		
		When $n$ is odd,  We will follow the same strategy as in previous theorem, i.e.,   we will show the stabilizer of null set $\emptyset$ equals to the group $S_n$.  Let $f \in Aut(H(n,k))$ such that $f(\emptyset)=\emptyset$. As $n-1$ is even hence $Aut(H(n-1,k))\cong H_{n-1}S_{n-1}$.  As $H(n-1,k)$ is induced subgraph of $H(n,k)$, let $f \vert_{[n-1]} \in S $, where $S$ is subgroup of $S_n$ and $S\cong S_{n-1}$. Let $f\vert_{[n-1]}=\sigma$.  Hence  we have $f(X)=\sigma(X)$ for all $X \subseteq [n-1]$, i.e., $|f(\{i\})|=1$ for all $i=1, \ldots , n-1$. Let $X \subset [n-1]$ and $|X|=k-1$. So $X \sim \{n\}$, hence $f(X)=\sigma(X) \sim f(\{n\})$ for every $(k-1)$-subset $X$ of $[n-1]$. Therefore $|f(\{n\})|=1$, i.e.,  $|f(\{i\})|=1$ for all $i=1, \ldots , n$.

		Let $f(\{i\})=\{ x_i \}$ for all $i \in \{1,2, \ldots , n\}$ and $\sigma \in S_n$ such that $\sigma(i)=x_i$. We will show that $f(X)=\sigma (X)$ for all $X \in \Omega_n$. Let $|X|=k+1$ and $i \in X$, hence $\{i\} \sim X$, i.e., $f(\{i\}) \sim f(X)$, i.e., $x_i \in f(X)$ for all $i \in X$, hence $f(X)=\sigma(X)$ for all $|X|=k+1$. Now let $|X|=k-1$ and  $i \notin X$, hence $\{i\} \sim X$, i.e., $f(\{i\}) \sim f(X)$, i.e., $x_i \notin f(X)$ for all $i \notin X$, hence $f(X)=\sigma(X)$ for all $|X|=k-1$.
		
		Let $|X|=2p+1$ such that $X$ has neighbours of size $k+1$. Consider the set $X'$ which contains $p+1$ many elements from $X$ and $k-p$ many elements from $X^c$. Therefore $|X'|=k+1$ and $X \sim X'$ and  $f(X')=\sigma(X')$. As $f(X)$ is common neighbour of all such $f(X')$, hence $f(X)=\sigma(X)$. Therefore $f(X)=\sigma(X)$ for all neighbours of of the sets of size $k+1$. Similarly we can prove for all the neighbours of the sets of size $k-1$. 
		
		Let  $f(X')=\sigma(X')$ for all $X' \in \Omega_n^{2p+1}$. Let $X \in \Omega_n$ such that $|X|=2m$ and $X$ has neighbour of size $2p+1$. Consider the set $X'$ which contains $p+m-\frac{k-1}{2}$ many elements from $X$ and $p-m+\frac{k+1}{2}$ many elements from $X^c$.  Therefore $|X'|=2p+1$ and $X \sim X'$ and  $f(X')=\sigma(X')$. As $f(X)$ is common neighbour of all such $f(X')$, hence $f(X)=\sigma(X)$. Therefore $f(X)=\sigma(X)$ for all neighbours of  the sets of size $2p+1$.
		
		Let  $f(X')=\sigma(X')$ for all $X' \in \Omega_n^{2m}$. Proceeding in similar manner we can prove $f(X)=\sigma(X)$ for all neighbours of the sets of size $2m$. Therefore through successive application of this method we have $f=\sigma \in S_n$. So the stabilizer of null set $\emptyset$ equals to the group $S_n$. By orbit-stabilizer theorem we have $|Aut(H(n,k))|=2^n n!$ and hence $Aut(H(n,k))=H_nS_n$. 
	\end{proof}
	
	When $k$ is even by the Theorem \ref{even k} we have $H(n,k)$ consists of two isomorphic connected components  $H'(n,k)$ and $H''(n,k)$ where  $H'(n,k)$ and $H''(n,k)$ are the induced subgraph of $H(n,k)$ by $\Omega'_n$ and $\Omega''_n$ respectively. Next, we will find the automorphism group of the subgraphs $H'(n,k)$ and $H''(n,k)$ when $n \neq 2k$. Let $H'_n=\{ \rho_X: X \in \Omega'_n \}$ and $H''_n=\{ \rho_X: X \in \Omega''_n \}$.

	{\theorem \label{Aut(H(n,k))3} $Aut(H'(n,k))=H'_n S_n$ and $Aut(H''(n,k))=H''_nS_n$, where $k$ is an even integer and $n \geq 2k+2$. }
	
	\begin{proof}
		We will prove this theorem for $H''(n,k)$ by orbit-stabilizer theorem.  Let $f \in Aut(H''(n,k))$ such that $f(\emptyset)=\emptyset$. Let $A_i$ denote the collection of sets which are at distance $i$ from $\emptyset$, hence $A_1=\Omega_n ^k$. So $f(A_i)=A_i$ for all $i$ and $|f(X)|=|X|$ for all $X \in A_1$.
		
		When $k=2$, $A_i$ is the set of size $2i$ for all $i$ in $H''(n,2)$. Similarly in $H'(n,2)$, distance between any singleton set and any set of size $(2p-1)$ is $p$. Therefore $|f(X)|=|X|$ for all $X \in \Omega_n$. Let $f(\{i\})=\{ x_i \}$ for all $i \in \{1,2, \ldots , n\}$ and $\sigma \in S_n$ such that $\sigma(i)=x_i$. Hence similarly as in the Theorem \ref{Aut(H(n,k))1} we can prove that $f \in S_n$. So the stabilizer of null set $\emptyset$ equals to the group $S_n$. By orbit-stabilizer theorem we have $|Aut(H''(n,k))|=2^{n-1} n!$ and hence $Aut(H''(n,k))=H''_nS_n$.

		\begin{figure}[ht]
			\centering
			\begin{center}
				\includegraphics[scale=1.5]{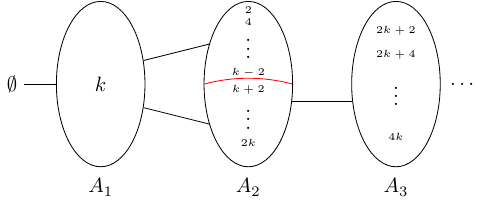}
				\label{pic}
			\end{center}
		\end{figure}

		Let $k\geq 4$. Note that as $n \geq 2k+2$, any set from $A_2$ of size $\geq k+2$ must have a neighbour in $A_3$, but the sets of size $\leq k-2$ has no neighbour in $A_3$. Therefore $|f(\{i,j\})| \leq k-2$.
		
		Let $X \in A_2$ such that $|X|=2p \leq k-2$. Number of neighbours of $X$ in $A_1$ is $\binom{2p}{p} \binom{n-2p}{k-p}=v_p$(say). As in the Theorem \ref{Aut(H(n,k))1} we have $\{v_p\}$ is monotone decreasing sequence for $1 \leq p \leq \frac{k-2}{2}$, $v_1 \gneq v_i $  for all $i=2, \ldots, \frac{k-2}{2}$ and $v_{\frac{k-2}{2}} \lneq v_i$ for all $i=1, \ldots ,\frac{k-4}{2}$. Therefore $|f(X)|=|X|$ for all $X \in A_2$ of size $2$ and $k-2$. 
		
		Let $X_1, X_2 \in A_2$ such that $|X_1|=|X_2|=2$.  The sets $\{a,b \}$ and $\{b,c\}$ have $\binom{n-3}{k-2}$ many and the sets $\{a,b\}$ and $\{c,d\}$ have $\binom{n-4}{k-2}$ many common neighbours of size $(k-2)$ for all different $a,b,c,d \in [n]$. Hence if $|X_1 \cap X_2|=1$ then $|f(X_1) \cap f(X_2)|=1$ and if $|X_1 \cap X_2|=0$ then $|f(X_1) \cap f(X_2)|=0$. Therefore we can take $f(\{ i,j\})=\{x_i,x_j\}$ and  $\sigma \in S_n$ such that $\sigma(i)=x_i$ for all $i=1,\ldots ,n$, i,e, $f(X)=\sigma(X)$ for all $|X|=2$. Now we can  proceed similarly as in Theorem \ref{Aut(H(n,k))1}, hence we have  $Aut(H''(n,k))=H''_nS_n$. As $H'(n,k)\cong H''(n,k)$, so $Aut(H'(n,k))=H'_nS_n$.
	\end{proof}

	{\theorem \label{Aut(H(n,k))4} $Aut(H'(n,k))=H'_n S_n$ and $Aut(H''(n,k))=H''_nS_n$, where $k$ is an even integer and $n \leq 2k-2$. }
	
	\begin{proof}
		We will prove this theorem by the method of induction on $n$ for the graph $H''(n,k)$. Let $ n=k+1$ and $f \in Aut(H''(k+1,k))$ such that $f(\emptyset)=\emptyset$. 
		Let $X \in \Omega''_n$. Therefore distance between $\emptyset$ and $X$ is $2p+1$ if $|X|=k-2p$ and is $2p$ if $|X|=2p$, where $ 0\leq p \leq \frac{k-2}{2}$. Therefore $|f(X)|=|X|$ for all $X \in \Omega''_n$.  Let $f(\{ i,j\})=\{x_i,x_j\}$ and  $\sigma \in S_n$ such that $\sigma(i)=x_i$ for all $i=1,\ldots ,n$, i,e, $f(X)=\sigma(X)$ for all $X \in \Omega_n^2$.
		Hence proceeding as Theorem \ref{Aut(H(n,k))2} we have $Aut(H''(k+1,k))=H''_{k+1}S_{k+1}$.
		
		Let the statement is true for $n-1$, now we will prove this for $n$. Let $f \in Aut(H''(n,k))$ such that $f(\emptyset)=\emptyset$. By induction hypothesis  $Aut(H(n-1,k))\cong H_{n-1}S_{n-1}$. As $H(n-1,k)$ is induced subgraph of $H(n,k)$, let $f \vert_{[n-1]} \in S $, where $S$ is subgroup of $S_n$ and $S\cong S_{n-1}$. Let $f\vert_{[n-1]}=\sigma$.  Hence  we have $f(X)=\sigma(X)$ for all $X \subseteq [n-1]$,  i.e., $|f(\{i,j\})|=2$ for all $i,j \in \{1, \ldots , n-1\}$. 
		
		Let $X \subset [n-1]$ such that $X=\{x_1,x_2, \ldots , x_k\}$. So $X \sim \{x_i,n\}$ for all $i=1,2, \ldots , k$, hence $f(X)=\sigma(X) \sim f(\{x_i, n\})$ for every $(k-1)$-subset $X$ of $[n-1]$. As $X$ is arbitrary $k$-subset of $[n-1]$,  hence $|f(\{i,n\})|=2$ for all $i=1,2, \ldots ,n-1$, i.e.,  $|f(\{i,j\})|=2$ for all $i,j=1, \ldots , n$.  Let $f(\{ i,j\})=\{x_i,x_j\}$ and  $\sigma \in S_n$ such that $\sigma(i)=x_i$ for all $i=1,\ldots ,n$, i,e, $f(X)=\sigma(X)$ for all $X \in \Omega_n^2$. Hence proceeding as Theorem \ref{Aut(H(n,k))2} we have $Aut(H''(n,k))=H''_{n}S_{n}$. As $H'(n,k)\cong H''(n,k)$, so $Aut(H'(n,k))=H'_nS_n$.

	\end{proof}

	\section{Geodesic Transitivity}	
	
	In this section we give full characterization of all $n$ and $k$ for which $H(n,k)$ is geodesic transitive.

	{\theorem \cite{Hu_Wang_Xu} For any odd integer $k$, $d(H(n,k))=\lceil\frac{n-1}{k}\rceil +1$ if $n \geq 2k-1 $;  $d(H(n,k))=\lceil \frac{n-1}{n-k} \rceil +1$ if $n \leq 2k-2$.}
	
	Note that when $n$ is even, $d(H(n,1))=n$, i.e,  $d(H(n,n-1))=n$. When $k$ is even by the Theorem \ref{even k} $H(n,k)$ contains two isomorphic components $H'(n,k)$ and $H''(n,k)$.  In $H''(n,2)$ the sets of size $n$ or $n-1$ are at maximum distance from $\emptyset$ and the distance is $\lfloor \frac{n}{2} \rfloor $ depends on the parity of $n$.  Hence diameter of each connected components of $H(n,2)$ is $\lfloor \frac{n}{2} \rfloor $. When $n-1$ is even,  in $H''(n,n-1)$ the sets of size $\frac{n-1}{2}$ or $\frac{n+1}{2}$ are at maximum distance from $\emptyset$ depends on parity of  $\frac{n-1}{2}$ and $\frac{n+1}{2}$,  and the distance is $\frac{n-1}{2}$. Hence diameter of each connected components of $H(n,n-1)$ is $\frac{n-1}{2}$.

	{\theorem $H(n,1)$ and $H(n,2)$ are  geodesic transitive.}
	
	\pf Consider two $l$-geodesics $G_1:\emptyset \sim \{x_1\} \sim \{x_1,x_2\}\sim \cdots \sim \{x_1,x_2,\ldots ,x_l\}$ and $G_2:\emptyset \sim \{y_1\} \sim \{y_1,y_2\}\sim \cdots \sim \{y_1,y_2,\ldots ,y_l\}$ where $1 \leq l \leq n$ in $H(n,1)$.  Consider the permutation $\sigma \in S_n$ such that $\sigma(x_i)=y_i$, for all $i=1,2,\ldots ,l$. Hence $\sigma(G_1)=G_2$. Hence $H(n,1)$ is geodesic transitive. 
	
	Arguing as above we have $H''(n,2)$ is  geodesic transitive. 	As $H(n,2)$ consists of two isomorphic connected components, hence $H(n,2)$ is  geodesic transitive. \qed

	{\theorem 
		$H(n,n-1)$ is geodesic transitive.}
	
	\pf If $n$ is even, i.e., $n-1$ is odd, then by the Theorem \ref{H(n,n-k)} we have $H(n,n-1) \cong H(n,1)$. Hence $H(n,n-1)$ is  geodesic transitive.

	Now consider the case when $n$ is odd and $\frac{n+1}{2}$ is even, i.e, $H(n,n-1)$ has two isomorphic components $H'(n,n-1)$ and $H''(n,n-1)$. Any $\frac{n-1}{2}$-geodesic  starting from $\emptyset$ must be of the form $\emptyset=A_1 \sim B_1 \sim A_2\sim \cdots  \sim A_{\frac{n+1}{4}} \sim B_{\frac{n+1}{4}} $ where, $|A_i|=2(i-1)$ and $|B_i|=(n-1)-2(i-1)$ for all $i$. From the adjacency condition we have the following properties:
	\begin{itemize}
		\item $A_1\subset A_2 \subset \cdots \subset A_{\frac{n+1}{4}}$.
		\item $B_1 \supset B_2\supset \cdots \supset B_{\frac{n+1}{4}}$.
		\item $|A_i \cap B_{i-1} |=1$ for all $i$.
		\item $A_i \cap B_i=\emptyset$ for all $i$.
	\end{itemize}
	Let $G_1$ and $G_2$ be two $i$-geodesics  such that: 
	$G_1: \emptyset \sim \{x_1, \ldots ,x_{n-2},x_{n-1}\} \sim \{x_{n-1},x_n\} \sim \{x_1,\ldots, x_{n-4},x_{n-3} \}\sim \{ x_{n-3},x_{n-2},x_{n-1},x_n \}\sim  \cdots  \sim \{x_1, \ldots ,x_{n-2i+1} \}$ and\\
	$G_2: \emptyset \sim \{y_1, \ldots ,y_{n-2},y_{n-1}\} \sim \{y_{n-1},y_n\} \sim \{y_1,\ldots, y_{n-4},y_{n-3} \}\sim \{ y_{n-3},y_{n-2},y_{n-1},y_n \}\sim  \cdots  \sim \{y_1, \ldots ,y_{n-2i+1} \}$. Consider the permutation $\sigma \in S_n$ such that $\sigma(x_i)=y_i$ for all $i=1,2,\ldots ,n$. Hence  $\sigma(G_1)=G_2$. Hence $H(n,n-1)$ is  geodesic transitive.
	
	\qed

	{\theorem $H(n,k)$ is not $2$-geodesic transitive for $k\geq 3$, $n\neq k+1$ and $(n,k) \neq (5,3),~(6,4),(7,4)$.} 
	
	\pf Let $n=2k$ and $X \in \Omega_n$. If $|X|=2$ or $X=[n]$ then number of neighbours of $X$ in $\Omega_n^k$ is $2\binom{2k-2}{k-1}=a$(say) or  $\binom{2k}{k}=b$(say) respectively. As $k \geq 3$, hence $a\neq b$. Hence for all $\varphi \in Aut(H(2k,k))$ with $\varphi(\emptyset)=\emptyset$ and $|X|=2$ imply $\varphi(X)\neq A$. Now consider the $2$-geodesics $G_1:\varphi \sim \{x_1,\ldots ,x_k\} \sim \{x_k,x_{k+1}\}$, and $G_2:\varphi \sim \{x_1,\ldots ,x_k\} \sim A$. There is no automorphism which maps $G_1$ to $G_2$, hence $H(2k,k)$ is not $2$ geodesic transitive.

	Let $n \neq 2k$, $k \neq 4$, i.e, $n \geq k+2$ and $(n,k)\neq (5,3)$. Consider the $2$-geodesics $G_1: \emptyset \sim \{x_1,\ldots ,x_k\}\sim \{x_k,x_{k+1}\}$ and  $G_3: \emptyset \sim \{x_1,\ldots ,x_k\}\sim \{x_{k-1},x_{k},x_{k+1},x_{k+2}\}$. From the above sections we have any set of size $2$ can be mapped to a set of size $n-1$. As $(n,k)\neq (5,3)$, i.e., $n-1 \neq 4$, hence there is no automorphism which maps $G_1$ to $G_3$, hence $H(n,k)$ is not $2$ geodesic transitive.
	
	Now let  $k=4$ and $n\geq 9$. Consider the $2$-geodesics $G_4: \emptyset \sim \{x_1,x_2,x_3 ,x_4\}\sim \{x_4,x_{5}\}$ and  $G_5: \emptyset \sim \{x_1,x_2,x_3 ,x_4\}\sim \{x_{2},x_{3},x_{4},x_{5},x_6,x_7\}$. As $n\geq 9$, i.e., $n-1 \neq 6$, hence there is no automorphism which maps $G_4$ to $G_5$, hence $H(n,4)$ is not $2$ geodesic transitive. \qed
	
	{\note When $(n,k)=(5,3),(6,4),(7,4)$, SAGE \cite{sagemath} computation shows that $H(n,k)$ is geodesic transitive.}

	\section{Open Issues}
	In this paper we have found the full automorphism group except for the case $n=2k$, so this case is still an open issue. M. Afkhami \textit{et.al.} \cite{New_Cayley} generalized $H(n, k) $ and defined a new family of Cayley graph, so the automorphism group of this family can be an interesting topic to investigate.
	Exploring the fixing number of this family of  graphs, a crucial parameter contingent upon the automorphism group, could be an interesting topic for additional research.
	
	\section{Acknowledgement}
	The author is supported by the Ph.D. Fellowship of CSIR (File No. 08/155(0086)/2020-
	EMR-I), Government of India. Special thanks to Angsuman Das and Anubrato Bhattacharyya for their valuable suggestions.

	\begin{center}
		\textbf{Appendix A\\
			\vspace{.2in}
			$\Gamma_{n+1}$ is a subgroup of $\Gamma$}
	\end{center}
	\pf We have $f_{\sigma}: \Omega_n \rightarrow \Omega_n$ such that
	\begin{equation}
	f_{\sigma }(X)=\left\lbrace \begin{array}{ll}
	\sigma(X) &  \mbox{~when~}  t \notin X\\
	
	(\sigma(X\setminus  \{t\} ))^c & \mbox{~when~}  t \in X.
	\end{array} \right. 
	\end{equation}
	where $t=\sigma^{-1}(n+1)$. 
	$\Gamma_{n+1}=\{ f_{\sigma} : \sigma \in S_{n+1} \}$. Let $f_{\sigma_1}, f_{\sigma_2} \in \Gamma_{n+1}$, we will prove that $f_{\sigma_1} \circ f_{\sigma_2}=f_{\sigma_1 \circ \sigma_2}$. Let $\sigma_1(t_1)=n+1$, $\sigma_2(t_2)=n+1$ and $\sigma_2(t_3)=t_1$, hence $(\sigma_1 \circ \sigma_2)(t_3)=n+1$. Let $X \in \Omega_n$. Note that complement of $X$ in $[n]$ is denoted by $X^c$ and complement of $X$ in $[n+1]$ is denoted by $[n+1]\setminus X$.
	
	\textbf{Case 1:} Let $t_3 \in X$. So $f_{\sigma_1 \circ \sigma_2}(X)=((\sigma_1 \circ \sigma_2)(X \setminus \{t_3\}))^c=((\sigma_1 \circ \sigma_2)(X)\setminus \{n+1\})^c.$
	
	
	Let $t_2 \in X$. Now $t_3 \in X$ implies $t_1=\sigma_2(t_3)\in \sigma_2(X)$, hence $t_1 \notin (\sigma_2(X\setminus \{t_2\}))^c$.
	So $(f_{\sigma_1} \circ f_{\sigma_2})(X)=f_{\sigma_1}((\sigma_2(X\setminus \{t_2\}))^c)=\sigma_1((\sigma_2(X\setminus \{t_2\}))^c)=((\sigma_1 \circ \sigma_2)(X)\setminus \{n+1\})^c.$
	
	
	\textbf{Case 2:} Let $t_3 \notin X$. This case can be done similarly.

	Therefore  $f^{\sigma_1} \circ f^{\sigma_2}=f^{\sigma_1 \circ \sigma_2}$ and  $(f_{\sigma})^{-1}=f_{\sigma^{-1}}$ and hence $\Gamma_{n+1}$ is a subgroup of $\Gamma$. \qed

	\begin{center}
		\textbf{Appendix B\\
			\vspace{.1in}
			$\Gamma^{n+1}$ is a subgroup of $\Gamma$}
	\end{center}
	\pf Let $\sigma \in S_{n+1}$ and $\sigma(t)=n+1$. If $t \in [n]$ then we have $f^{\sigma}: \Omega_n \rightarrow \Omega_n$ such that
	
	\begin{equation}
	f^{\sigma }(X)=\left\lbrace \begin{array}{ll}
	\sigma(X) &  \mbox{~if~} |X| \mbox{~is even  and~}t \notin X\\
	(\sigma(X \setminus \{t\}))^c & \mbox{~if~} |X| \mbox{~is even  and~}t \in X\\
	(\sigma(X \cup \{n+1\}))^c & \mbox{~if~} |X| \mbox{~is odd  and~}t \notin X\\
	\sigma(X\cup \{n+1\} \setminus \{t\}) &  \mbox{~if~} |X| \mbox{~is odd  and~}t \in X.
	\end{array} \right. 
	\end{equation} 
	
	If $\sigma(n+1)=n+1$ then we have $f^{\sigma}(X)=\sigma(X)$. $\Gamma^{n+1}=\{ f^{\sigma} : \sigma \in S_{n+1} \}$. Let $f^{\sigma_1}, f^{\sigma_2} \in \Gamma^{n+1}$, we will prove that $f^{\sigma_1} \circ f^{\sigma_2}=f^{\sigma_1 \circ \sigma_2}$.  Let $\sigma_1(t_1)=n+1$, $\sigma_2(t_2)=n+1$ and $\sigma_2(t_3)=t_1$.  Hence $(\sigma_1 \circ \sigma_2)(t_3)=n+1$. Note that complement of $X$ in $[n]$ is denoted by $X^c$ and complement of $X$ in $[n+1]$ is denoted by $[n+1]\setminus X$. If both $t_1=t_2=n+1$ then the result is obvious.

	\textbf{Case 1:} Let $t_1=n+1$ and  $t_2 \in [n]$ , i.e., $\sigma_1(n+1)=n+1$, $t_2=t_3$. 
	\begin{itemize}
		\item Let $X\in \Omega''_n$ and $t_2\in X$.
		$f^{\sigma_1 \circ \sigma_2}(X)=( (\sigma_1 \circ \sigma_2)(X \setminus \{t_2\} )^c$.
		$(f^{\sigma_1} \circ f^{\sigma_2})(X)=f^{\sigma_1}((\sigma_2(X \setminus \{t_2\} ))^c)=\sigma_1( (\sigma_2(X \setminus \{t_2\} ))^c )=( (\sigma_1 \circ \sigma_2)(X \setminus \{t_2\} )^c$.
		
		\item  Let $X\in \Omega''_n$ and $t_2\notin X$. 
		$f^{\sigma_1 \circ \sigma_2}(X)=(\sigma_1 \circ \sigma_2)(X )=(f^{\sigma_1} \circ f^{\sigma_2})(X)$.
		
		\item  Let $X\in \Omega'_n$ and $t_2\in X$. 
		$f^{\sigma_1 \circ \sigma_2}(X)=(\sigma_1 \circ \sigma_2)(X \cup \{n+1\} \setminus \{t_2\})$.
		$(f^{\sigma_1} \circ f^{\sigma_2})(X)=f^{\sigma_1}(\sigma_2(X \cup \{n+1\} \setminus \{t_2\}))=\sigma_1( \sigma_2(X \cup \{n+1\} \setminus \{t_2\}) )=(\sigma_1 \circ \sigma_2)(X \cup \{n+1\} \setminus \{t_2\})$.
		
		\item  Let $X\in \Omega'_n$ and $t_2\notin X$. 
		$f^{\sigma_1 \circ \sigma_2}(X)=((\sigma_1 \circ \sigma_2)(X \cup \{n+1\}))^c$.
		$(f^{\sigma_1} \circ f^{\sigma_2})(X)=f^{\sigma_1}((\sigma_2 ( X \cup \{n+1\} ))^c)=\sigma_1( (\sigma_2 ( X \cup \{n+1\} ))^c )
		=((\sigma_1 \circ \sigma_2)(X \cup \{n+1\}))^c$.
		
	\end{itemize}

	\textbf{Case 2:} Let $t_1 \in [n]$ and $t_2=n+1$. This case can be done as previous.

	\textbf{Case 3:} Let  $t_1,t_2 \in [n]$.

	\textbf{Subcase 1:} Let $t_3=n+1$, i.e., $(\sigma_1 \circ \sigma_2)(n+1)=n+1$. Hence $f^{\sigma_1 \circ \sigma_2}(X)=(\sigma_1 \circ \sigma_2)(X)$.  As $t_3 \notin X$, hence $t_1 \notin \sigma_2(X)$ and $t_1 \in \sigma_2(X \cup \{n+1\}) $. 
	
	\begin{itemize}

		\item Let $X\in \Omega''_n$ and $t_2\in X$. $(f^{\sigma_1} \circ f^{\sigma_2})(X)=f^{\sigma_1}( (\sigma_2(X \setminus \{t_2\})^c )=\sigma_1( ((\sigma_2(X \setminus \{t_2\})^c\setminus \{t_1\})^c )=(\sigma_1 \circ \sigma_2)(X)$.
		
		\item Let $X\in \Omega''_n$ and $t_2\notin X$. $(f^{\sigma_1} \circ f^{\sigma_2})(X)=f^{\sigma_1}(\sigma_2(X))=\sigma_1(\sigma_2(X) )=(\sigma_1 \circ \sigma_2)(X)$.

		\item Let $X\in \Omega'_n$ and $t_2\in X$. $(f^{\sigma_1} \circ f^{\sigma_2})(X)=f^{\sigma_1}( \sigma_2(X \cup \{n+1\} \setminus \{t_2\}) )=f^{\sigma_1}(\sigma_2(X)\cup \{t_1\} \setminus \{n+1\})
		=(\sigma_1 \circ \sigma_2)(X).$

		\item Let $X\in \Omega'_n$ and $t_2\notin X$. $(f^{\sigma_1} \circ f^{\sigma_2})(X)=f^{\sigma_1}( (\sigma_2(X \cup \{n+1\}))^c )=(  (\sigma_2(X \cup \{n+1\}))^c \cup \{n+1\} )^c=(\sigma_1 \circ \sigma_2)(X).$
	\end{itemize}
	\textbf{Subcase 2:} Let $t_3 \in [n].$ Then we have two cases, either $t_3 \in X$ or $t_3 \notin X$. This case can be done as previous.
	
	Therefore  $f^{\sigma_1} \circ f^{\sigma_2}=f^{\sigma_1 \circ \sigma_2}$ and $(f^{\sigma})^{-1}=f^{\sigma^{-1}}$ and hence $\Gamma^{n+1}$ is a subgroup of $\Gamma$. \qed

	\begin{center}
		\textbf{Appendix C}
	\end{center}
	\vspace{.2in}
	The sequence $\{u_p\}$ where $u_p=\binom{2p}{p} \binom{n-2p}{k-p}$, is strictly decreasing for $p=1,2, \ldots , \frac{k-3}{2}$ and $n \geq 2k+2$.
	
	\pf $\frac{u_p}{u_{p+1}}=\frac{\binom{2p}{p} \binom{n-2p}{k-p}}{\binom{2p+2}{p+1} \binom{n-2p-2}{k-p-1}}=\frac{(p+1)(n-2p)(n-2p-1)}{2(2p+1)(k-p)(n-k-p)}$. As $n \geq 2k+2$, so $n-2k>0$, i.e., $n^2+4p^2>4nk$, i.e., $n^2p+4k^2p>4nkp $ and $2nk>n$. 
	
	Therefore we have, $n^2+4k^2p+2nk+(n-2p)^2+np+2p+2k^2> 4nkp+n$, i.e., $(p+1)\{(n-2p)^2-(n-2p)\}>2(2p+1)(nk-np-k^2+p^2)$, hence $u_p>u_{p+1}$. \qed

\end{document}